\documentclass[a4paper,11pt]{article}
\usepackage{amsmath}
\usepackage{amssymb}
\usepackage{mathtools}
\usepackage{theorem}
\usepackage{euscript}
\usepackage{epic,eepic}
\usepackage{graphicx}
\usepackage{xcolor}
\usepackage{hyperref}
\hypersetup
{
    colorlinks,
    citecolor=green,
    filecolor=black,
    linkcolor=blue,
    urlcolor=blue
}
\hypersetup{linktocpage}

\PassOptionsToPackage{normalem}{ulem}
\usepackage{amsmath}
\usepackage{amssymb}
\usepackage{mathtools}
\usepackage{theorem}
\usepackage{euscript}
\usepackage{epic,eepic}
\usepackage{graphicx}
\usepackage{setspace}
\usepackage{comment}
\PassOptionsToPackage{normalem}{ulem}

\newcommand{\menge}[2]{\big\{{#1} \;|\; {#2}\big\}}

\newcommand{\emp}{\ensuremath{{\varnothing}}}

\newcommand{\scal}[2]{\left\langle{#1}\mid {#2} \right\rangle}

\newcommand{\vuo}{\ensuremath{\mbox{\footnotesize$\square$}}}
\newcommand{\vuoO}{\ensuremath{\mbox{\footnotesize$\square_O$}}}

\newcommand{\HH}{\ensuremath{\mathcal H}}
\newcommand{\GG}{\ensuremath{\mathcal G}}
\newcommand{\YY}{\ensuremath{\mathcal Y}}
\newcommand{\XX}{\ensuremath{\mathcal X}}

\newcommand{\KKK}{\ensuremath{\boldsymbol{\mathcal K}}}

\newcommand{\MM}{\ensuremath{\boldsymbol{M}}}

\newcommand{\sri}{\ensuremath{\operatorname{sri}}}

\newcommand{\RR}{\ensuremath{\mathbb R}}
\newcommand{\KK}{\ensuremath{\mathcal K}}

\newcommand{\NN}{\ensuremath{\mathbb N}}

\newcommand{\dom}{\ensuremath{\operatorname{dom}}}

\newcommand{\prox}{\ensuremath{\operatorname{prox}}}

\newcommand{\ran}{\ensuremath{\operatorname{ran}}}
\newcommand{\zer}{\ensuremath{\operatorname{zer}}}
\newcommand{\gra}{\ensuremath{\operatorname{gra}}}

\newcommand{\vv}{\ensuremath{\boldsymbol{v}}}
\newcommand{\ww}{\ensuremath{\boldsymbol{w}}}

\newcommand{\xx}{\ensuremath{\mathsf{x}}}
\newcommand{\af}{\ensuremath{\mathsf{a}}}

\newcommand{\yy}{\ensuremath{\mathsf{y}}}

\newcommand{\uu}{\ensuremath{\boldsymbol{u}}}

\newcommand{\LL}{\ensuremath{\boldsymbol{L}}}

\newcommand{\SSS}{\ensuremath{\boldsymbol{S}}}
\newcommand{\QQ}{\ensuremath{\boldsymbol{Q}}}

\newcommand{\Id}{\ensuremath{\operatorname{Id}}}

\newcommand{\weakly}{\ensuremath{\rightharpoonup}}

\newtheorem{theorem}{Theorem}[section]
\newtheorem{lemma}[theorem]{Lemma}

\newtheorem{definition}[theorem]{Definition}
\theoremstyle{plain}{\theorembodyfont{\rmfamily}
}
\theoremstyle{plain}{\theorembodyfont{\rmfamily}
}
\theoremstyle{plain}{\theorembodyfont{\rmfamily}
\newtheorem{algorithm}[theorem]{Algorithm}}
\theoremstyle{plain}{\theorembodyfont{\rmfamily}
}
\theoremstyle{plain}{\theorembodyfont{\rmfamily}
\newtheorem{problem}[theorem]{Problem}}
\theoremstyle{plain}{\theorembodyfont{\rmfamily}
\newtheorem{remark}[theorem]{Remark}}
\theoremstyle{plain}{\theorembodyfont{\rmfamily}
\newtheorem{notation}[theorem]{Notation}}

\definecolor{labelkey}{rgb}{0,0.08,0.45}
\definecolor{refkey}{rgb}{0,0.6,0.0}
\definecolor{Brown}{rgb}{0.45,0.0,0.05}
\definecolor{dgreen}{rgb}{0.00,0.49,0.00}
\definecolor{dblue}{rgb}{0,0.08,0.75}
\numberwithin{equation}{section}

\title{An algorithm for solving structured problems involving $L$-infimal convolutions}

\author{ V\~u C\^ong B$\grave{\text{\u{a}}}$ng$^1$ and V\~u Xu\^an  Nh\^am$^2$   \\[5mm] 
$^1$ \small Belgium Research Center (BeRC) - Huawei, Leuven, Belgium\\
 $^2$ \small University of Transport Technology, Hanoi, Vietnam  \\       
\small bangcvvn@gmail.com; nhamvx@utt.edu.vn}  
\date{}
\begin{document}
\maketitle
\begin{abstract} 
We develop an operator-splitting framework for solving structured minimization problems involving linear compositions of $L$-infimal convolutions. We introduce a primal-dual shadow Douglas–Rachford splitting method and establish its convergence. As a result, we derive a fully primal-dual splitting algorithm for the class of the aforementioned  minimization problems.
\end{abstract}


\section{Introduction}
 Over the past two decades, a wide range of algorithms has been developed in the literature for solving monotone inclusion problems, ranging from purely primal methods to primal-dual splitting approaches.
The primal dual splitting algorithm for solving monotone inclusions was first introduced in \cite{siop2} based on a reformulation of the underlying inclusions as a monotone-skew model. In the later years, this primal dual monotone inclusion was extended extensively in the literature, see \cite{plc6,vu2013,Bot2013,acta2,Dao26,Buicomb,Malistky26} and the references therein, for instances. Monotone inclusions provide a unifying mathematical framework for a broad class of problems in applied mathematics, including optimization, variational analysis, equilibrium theory, inverse problems, transport optimization, and machine learning. A typical application is to minimize the structured convex optimization problem involving the infimal convolutions; see \cite{plc6, vu2013,Buicomb,Becker,MinYan, Fadili} for examples. Recall that The infimal convolution of the two functions $\ell$ and $g$ from a Hilbert space $\GG$ to $\left]-\infty,+\infty\right]$ is \cite{livre1}
 \begin{equation}
 \ell\;\vuo\; g\colon x \mapsto \inf_{y\in\GG}(\ell(y)+g(x-y)).
\end{equation}
 This definition requires $\ell$ and $g$ define on the same space $\GG$. Now, let $f$ be a function from  $\HH$ to $\left]-\infty,+\infty\right]$ and $L\colon \HH \to \GG$ be an operator,  define  $L$-infimal convolution of $f$ and $g$ as
 \begin{equation}
 f\;\vuo_L\; g\colon \GG\ni x \mapsto \inf_{y\in\HH}(f(y)+g(x-Ly)).
\end{equation}
The $L$-infimal convolution operation $\vuo_L$ is not symmetric and it is an extension of the infimal convolution. In this work, we focus on the following problem which extends primal problem in \cite[Problem 4.1]{plc6}.
\begin{problem}
\label{prob:app0}
Let $\HH$ be a real Hilbert space and $z\in\HH$, let $f\in\Gamma_0(\HH)$, and let $h\colon\HH\to\RR$ be convex and
differentiable with a $\beta^{-1}$-Lipschitzian gradient for 
some $\beta\in\left]0,+\infty\right[$.
Let $m$ be a strictly positive integer.
For every $i \in \{1,\ldots, m\}$, let $\GG_i$, $\mathcal{X}_i$ and $\mathcal{Y}_i$ be a real Hilbert spaces, let $W_i\colon \mathcal{X}_i\to \GG_i$ and $O_i\colon\mathcal{Y}_i\to\GG_i$ be nonzero bounded linear operators,
and $r_i\in\GG_i$, 
 let $g_i \in \Gamma_0(\YY_i)$, 
let $\ell_i \in \Gamma_0(\XX_i)$ be $\mu_i$-strongly convex, for some
$\mu_i\in\left]0,+\infty\right[$ such that  $0\in \sri(\ran O_{i}^*- \dom g^{*}_i)$.
For every $i \in \{1,\ldots, m\}$, let  $L_i\colon\HH\to\GG_i$ is a nonzero bounded linear operator, and define
\begin{equation}
     W_i \triangleright \ell_i\colon w\mapsto \inf_{\substack{u \in \XX_i \\ W_iu = w}} \ell_i(u),
    \; 
\end{equation}
which is the infimal postcomposition of $\ell_i$ by $W_i$.
The problem is to
\begin{equation}\label{e:pri10}
\underset{x\in\HH}{\text{minimize}} \;
f(x)+\sum_{i=1}^m\big( g_i\;\vuoO_i\;(W_i \triangleright\ell_i)\big)(L_ix-r_i) + h(x) -\scal{z}{x}.
\end{equation}
\end{problem}
The dual problem of \eqref{e:pri10} will be stated in the last section. The problem \eqref{e:pri10} models various structured minimization problem in the literature. 
The objective of this paper is to develop the operator splitting 
technique the problem \eqref{e:pri10} and its dual problem. To  this ends, we first revisit an instance of the ones in \cite{Buicomb,acta2}. 
\begin{problem}\label{p:Prob0}
Let $m$ be a strictly positive integer, 
let $\beta$, $\beta_0$ and $(\mu_i)_{1\leq i\leq m}$ be in $\left]0,+\infty\right[$.
Let $C\colon\HH\to\HH$ be  $\beta$-cocoercive, let $Q\colon\HH\to\HH$  $\beta_0$-Lipschitzian monotone,
let  $A\colon \HH\to 2^{\HH}$ be maximally monotone, and let $z\in\HH$.
Let $(\GG_i)_{1\leq i\leq m}$ be real Hilbert spaces. 
For every $i\in\{1,\ldots,m\}$, let $B_i\colon\GG_i\to 2^{\GG_i}$ be  maximally monotone, 
 let $D_i\colon\GG_i\to 2^{\GG_i}$ be such that $D_{i}^{-1}$  is $\mu_i$-cocoercive,  let $r_i\in\GG_i$, and
 let $L_i\colon\HH\to \GG_i$ be a bounded linear operator such that $0\not=\sum_{i=1}^m\|L_i\|^2$. 
Suppose that 
\begin{equation}
\label{c:1}
z \in \ran\bigg(A+ \sum_{i=1}^m L_{i}^* (B_{i}\vuo  D_i)( L_i\cdot-r_i) + C +Q\bigg).
\end{equation}
The primal inclusion is to find $\overline{x}$ such that
\begin{equation}
\label{e:pri}
z \in A\overline{x}+ \sum_{i=1}^m L_{i}^* \big( B_i\vuo D_i\big)  (L_i \overline{x}-r_i)+ (C+Q)\overline{x},
\end{equation}
and the dual inclusion is to find $(\overline{v}_i)_{1\leq i\leq m} \in (\GG_i)_{1\leq i\leq m}$ such that 
\begin{equation}
\label{e:dual}
\big(\forall i\in\{1,\ldots,m\}\big) \; -r_i \in -L_i(A+C+Q)^{-1}\bigg(z-\sum_{i=1}^m L_{i}^*\overline{v}_i\bigg) + B_{i}^{-1} \overline{v}_i 
+D_{i}^{-1} \overline{v}_i.
\end{equation}
We denote by $\mathcal{P}$ and $\mathcal{D}$ the set of solutions to \eqref{e:pri} and \eqref{e:dual}, respectively. 
\end{problem}
Section \ref{s:NB} is devoted to notation and background. Section 
\ref{s:Algo} presents the algorithm for solving Problem \ref{p:Prob0}, we prove its convergences.  We present the dual problem of \eqref{e:pri10} as well as the algorithm for solving it  and its dual problem.
\section{Notations and Background}
\label{s:NB}
In this work, $\HH$, $(\GG_i)_{1\leq i\leq m}$,$(\XX_i)_{1\leq i\leq m}$, $(\YY_i)_{1\leq i\leq m}$ are real Hilbert spaces. Their scalar products and associcated norms are respectively denoted by 
$\scal{\cdot}{\cdot}$ and $\|\cdot\|$.
Let $A\colon\HH\to 2^{\HH}$ be a set-valued operator.
The domain of $A$ is denoted by $\dom(A)$ that 
is a set of all $x\in\HH$ such that $Ax\not= \emp$.
The range of $A$ is $\ran(A) = \menge{u\in\HH}{(\exists x\in\HH)\ u\in Ax }$. 
The graph of $A$ is 
$\gra(A) = \menge{(x,u)\in\HH\times\HH}{u\in Ax}$.  
The inverse of $A$ 
is $A^{-1}\colon u \mapsto \menge{x}{ u\in Ax}$.  
The zero set of $A$ is $\zer(A) = A^{-1}0$.
\begin{definition} Let $A\colon\HH\to 2^{\HH}$ be a set valued-operator.
\begin{enumerate}
\item $A$ is  monotone if 
\begin{equation}\label{oioi}
\big(\forall (x,u)\in \gra A\big)\big(\forall (y,v)\in\gra A\big)
\quad\scal{x-y}{u-v}\geq 0.
\end{equation}
\item $A$ is maximally monotone if it is monotone and  there exists no monotone operator $B$ such that $\gra(B)$ properly contains $\gra(A)$.
\item   $A$  is $\phi_A$-uniformly monotone, at $x\in\dom (A)$, if there 
 exists an increasing function $\phi_A\colon\left[0,\infty\right[\to \left[0,\infty\right]$ that vanishes only at $0$ such that 
 \begin{equation}\label{unif}
\big(\forall u\in Ax\big)\big(\forall (y,v)\in\gra A\big)
\quad\scal{x-y}{u-v}\geq \phi_A(\|y-x\|).
\end{equation}
If $\phi_A = \nu_A|\cdot|^2$  for some $\nu_A\in \left]0,\infty\right[$, then we say that $A$ is $\nu_A$-strongly monotone.
\item The resolvent of $A$ is 
\begin{equation}
 J_A=(\Id + A)^{-1},\;  
\end{equation}
where $\Id$ denotes the identity operator on $\HH$.  
\end{enumerate}
\end{definition}

\begin{notation}\label{No:0} Throughout, 
we define the Hilbert spaces $\GG=\bigoplus_{i=1}^m\GG_i$ and $\KK= \HH\bigoplus\GG$. The scalar products and the norm  on $\GG$ is 
\begin{equation}
    \scal{\cdot}{\cdot}_{\GG}\colon (\vv,\ww)\mapsto \sum_{i=1}^m\scal{v_i}{w_i}\;
    \text{and}\;
    \|\cdot\|_{\GG}\colon \vv \mapsto \sqrt{\sum_{i=1}^m\|v_i\|^2};
\end{equation}
and on $\KK$ is 
    \begin{equation}
    \scal{\cdot}{\cdot}_{\KK}\colon ((x,\vv),(y,\ww))\mapsto \scal{x}{y}+\scal{\vv}{\ww}_\GG\;
    \text{and}\;
    \|\cdot\|_{\KK}\colon (x,\vv) \mapsto \sqrt{\|x\|^2+\|\vv\|^{2}_\GG}.
\end{equation}
Let $\tau$ and $(\sigma_i)_{1\leq i\leq m}$ be strictly positive real numbers, we define the weighed norm  on $\KK$ by
\begin{equation}\label{e:weigh}
    \||\cdot\||\colon (x,\vv)\mapsto \sqrt{\frac{1}{2\tau} \|x\|^2 +\sum_{i=1}^m\frac{1}{2\sigma_i} \|v_i\|^2 }.
\end{equation}
\end{notation}
\begin{notation}\label{No:1} Define 
 \begin{equation}
 \begin{cases}
\LL\colon\HH\to\GG\colon x\mapsto (L_1x,\cdots,L_mx)\\
\SSS \colon \KKK\to\KKK \colon (x,v_1,\ldots, v_m) \mapsto  (Qx+\sum_{i=1}^m L_{i}^*v_i,  -L_1x, \ldots,  -L_m x ),\\
 \QQ \colon \KKK\to\KKK \colon (x,v_1,\ldots, v_m) \mapsto (Cx, D_{1}^{-1}v_1,\ldots,  D_{m}^{-1}v_m),\\
 \MM\colon \KKK\to 2^{\KKK}\colon (x,v_1,\ldots, v_m) \mapsto (-z+Ax, r_1+ B^{-1}_{1}v_1, \ldots, r_m+B^{-1}_mv_m).\notag
 \end{cases}
 \end{equation}
\end{notation}
\begin{lemma}{\rm \cite[Eqs.(3.12),(3.13),(3.21),(3.22)]{plc6}}\label{l:1} Under the condition \eqref{c:1}, $\zer(\MM+\SSS+\QQ)\not=\emp$. Moreover,
    \begin{equation}
        (\overline{x},\overline{v})\in \zer(\MM+\SSS+\QQ)
        \Longrightarrow \overline{x}\in\mathcal{P} \; \text{and}\; \overline{v}\in\mathcal{D}.
    \end{equation}
\end{lemma}
\begin{lemma}{\rm \cite[Lemma 2]{Polyak1987}}\label{lem:polyak} Let $(c_k)_{k\in\NN}$, $(\xi_k)_{k\in\NN}$, $(\chi_k)_{k\in\NN}$ and $(\eta_k)_{k\in\NN}$ be nonnegative sequences such that 
\begin{equation}
    \begin{cases}
       &(\forall k\in\NN)\; \xi_{k+1} \leq (1+\eta_k) \xi_k+c_k-\chi_k,\\
       &\sum_{k\in\NN} \eta_k <+\infty, \sum_{k\in\NN}c_k <+\infty.
    \end{cases}
\end{equation}
Then $(\xi_k)_{k\in\NN}$ converges to $\overline{\xi}\in \RR_+$ and $\sum_{k\in\NN}\chi_k < +\infty $.    
\end{lemma}
\section{Algorithm and convergences}
\label{s:Algo}
We propose the following algorithm for solving Problem \ref{p:Prob0}. The algorithm is full-splitting and it  evaluates the set-valued operators via their resolvents and the singled-valued operator via their values. 
\begin{algorithm} Let $\tau$ and $(\sigma_i)_{1\leq i\leq m}$ be strictly positive real numbers. Let $x_0\in\HH$ and, for each $i\in \{1,\ldots,m\}$, let $v_{i,0}$ in $\GG_i$. Iterate
\begin{equation}
\label{algo:primal-dual}
\begin{array}{l}
\left\lfloor
\begin{array}{l}
x_{k+1} = J_{\tau  A} \Big(x_k-\tau \big(\sum_{i=1}^mL_{i}^*v_{i,k}+ Cx_k+2Qx_k-Qx_{k-1} -z\big)\Big)\\
\operatorname{For \; \text{$i =1,\ldots,m$}}\\
\begin{array}{l}
\left\lfloor
\begin{array}{l}
v_{i,k+1} = J_{\sigma_i  B^{-1}_{i}} \Big(v_{i,k}+\sigma_i\big(L_ix_{k+1}-D_{i}^{-1}v_{i,k} - r_i\big)\Big) + \sigma_i(L_ix_{k+1}-L_ix_k).\\
\end{array}
\right.\\[2mm]
\end{array}\\
\end{array}
\right.\\[2mm]
\end{array}
\end{equation}
\end{algorithm}

The main convergence results of the algorithm \eqref{algo:primal-dual} can be now stated. We prove the weak convergence of the primal dual sequence as well as the strong convergence under the uniform monotonicity condition. 
\begin{theorem}\label{t:main1}
Let $(x_k)_{k\in\NN}$ and $(v_{1,k},\ldots,v_{m,k})_{k\in\NN}$ be primal-dual sequences generated by Algorithm \ref{algo:primal-dual}.
Suppose that the following condition is satisfied
\begin{equation}\label{e:constep}
    \begin{cases}
     \big(\forall i\in\{1,\ldots,m\}\big)\;   \sigma_i < \mu_i,\\
        \frac{1}{2\tau}-\frac{1}{4\beta}
 -\beta_0 - \sum_{i=1}^m\sigma_i\|L_i\|^2\Big( \frac{\sigma_i}{2\mu_i}+\frac{2\mu_i}{\mu_i-\sigma_i}-\frac{3}{2}\Big) > 0.
    \end{cases}
\end{equation}
Then the following holds for some $\overline{x}\in\mathcal{P}$ and $(\overline{v}_1,\ldots,\overline{v}_m)\in\mathcal{D}$.
\begin{enumerate}
\item \label{t:main1i} $\sum_{k\in\NN}\|x_{k+1}-x_k\|^2 < +\infty$ and $(\forall i\in\{1,\ldots,m\})\; \sum_{k\in\NN}\|v_{i,k+1}-v_{i,k}\|^2<+\infty$.
\item \label{t:maini2} $\sum_{k\in\NN}\|Cx_{k}-C\overline{x}\|^2 < +\infty$ and $(\forall i\in\{1,\ldots,m\})\; \sum_{k\in\NN}\|D_{i}^{-1}v_{i,k}-
D_{i}^{-1}\overline{v}_{i}\|^2<+\infty$.
    \item\label{t:main1ii} The sequence $(x_k)_{k\in\NN}$ converges weakly to $\overline{x}$. 
    \item \label{t:main1iii}The sequence $(v_{1,k},\ldots,v_{m,k})_{k\in\NN}$ converges weakly to $(\overline{v}_1,\ldots,\overline{v}_m)$.
    \item\label{t:mainiv} Suppose that $C$ or $A$ is uniformly monotone at $\overline{x}$, then $x_k\to \overline{x}$.
    \item \label{t:mainv} Suppose that $D_{j}^{-1}$ or $B_{j}^{-1}$ is uniformly at $\overline{v}_{j}$, for some $j\in\{1,\ldots,m\}$, then $v_{j,k}\to \overline{v}_j$.
\end{enumerate}
\end{theorem}
\begin{proof} We first observe that under condition 
\eqref{e:constep}, there exist $\beta'$ and $(\mu'_{i})_{1\leq i\leq m}$ such that 
\begin{equation}\label{e:constep1}
    \begin{cases}
     \big(\forall i\in\{1,\ldots,m\}\big)\;   \sigma_i < \mu'_i< \mu_i, 0< \beta' < \beta\\
        \frac{1}{2\tau}-\frac{1}{4\beta'}-\beta_0
  - \sum_{i=1}^m\sigma_i\|L_i\|^2\Big( \frac{\sigma_i}{2\mu'_i}+\frac{2\mu'_i}{\mu'_i-\sigma_i}-\frac{3}{2}\Big) > 0.
    \end{cases}
\end{equation}
Define the following set
\begin{equation}
    \overline{Z}=\menge{(x,v_1-\sigma_1 L_1x,\ldots, v_m-\sigma_m L_mx)}{(x,v_1,\ldots,v_m)\in Z}.
\end{equation}
Let $(x,\ww)\in\overline{Z}$. Then there exists
$(x,\vv) \in Z = \zer(\MM+\SSS+\QQ)$ such that $(\forall i\in\{1,\ldots,m\})\; w_i= v_i-\sigma_i L_ix$. Then, we have
\begin{equation}\label{e:split1}
    z-\sum_{i=1}^m L_{i}^*v_i \in Ax+Cx +Qx\;
    \text{and}
    \;
   (\forall i\in\{1,\ldots,m\})\; L_ix -D_{i}^{-1}v_i \in r_i + B_{i}^{-1}v_i
\end{equation}
Let $k\in\NN$.
We define $\xx_k= (x_k,\vv_k)$ where $\vv_k=(v_{1,k},\ldots,v_{m,k})$.
 By the definition of the resolvent, the primal update in Algorithm \ref{algo:primal-dual} yields 
    \begin{equation}
\label{e:priup1}
        \frac{x_k-x_{k+1}}{\tau} - Cx_k-(2Qx_k-Qx_{k-1})-\sum_{i=1}^m L_{i}^*v_{i,k} \in -z+ Ax_{k+1}.
    \end{equation} 
In view of $\LL$ in Notation \ref{No:1}, we observe that $\LL^*\vv_k=\sum_{i=1}^m L_{i}^*v_{i,k}$. Hence, \eqref{e:priup1} and the first inclusion in \eqref{e:split1} can be rewritten as 
    \begin{equation}
\label{e:pri1a}
\begin{cases}
        \frac{x_k-x_{k+1}}{\tau} - Cx_k-(2Qx_k-Qx_{k-1})-\LL^*\vv_{k} \in -z+ Ax_{k+1}\\
         -\LL^*\vv -Cx-Qx \in -z+ Ax.
         \end{cases}
    \end{equation}
Since $A$ is monotone, it follows from \eqref{e:pri1a} that 
    \begin{align}
   0&\leq \scal{ \frac{x_k-x_{k+1}}{\tau} - Cx_k-(2Qx_k-Qx_{k-1})-\LL^*\vv_k +\LL^*\vv+Cx+Qx}{x_{k+1}-x}   \notag\\
   &= \frac{1}{\tau}\scal{x_k-x_{k+1}}{x_{k+1}-x}+\scal{Cx-Cx_k}{x_{k+1}-x} +\scal{\LL^*\vv-\LL^*\vv_k}{x_{k+1}-x}\notag\\
   &\quad +\scal{Qx-Qx_k}{x_{k+1}-x}-\scal{Qx_k-Qx_{k-1}}{x_{k+1}-x}.
   \label{e:fire1a}
\end{align}
 Since $C$ is $\beta$-cocoercive and $Q$ is monotone, we have 
    \begin{equation}
    \begin{cases}
        0\leq \scal{x_{k+1}-x}{Cx_{k+1}-Cx} -\beta\|Cx_{k+1}-Cx\|^2\\
          0\leq \scal{x_{k+1}-x}{Qx_{k+1}-Qx}.
        \label{e:coco1}
        \end{cases}
    \end{equation}
    Moreover, we also have
    \begin{align}
        \scal{Qx_k-Qx_{k-1}}{x_{k+1}-x}&= \scal{Qx_k-Qx_{k-1}}{x_{k+1}-x_k}+\scal{Qx_k-Qx_{k-1}}{x_k-x}\notag\\
        &\geq -\beta_0\|x_k-x_{k-1}\|\|x_k-x_{k+1}\|+\scal{Qx_k-Qx_{k-1}}{x_k-x}\notag\\
        &-\frac{\beta_0}{2}\big(\|x_k-x_{k-1}\|^2+\|x_k-x_{k+1}\|^2 \big)+\scal{Qx_k-Qx_{k-1}}{x_k-x}.
        \label{e:Lipds1}
    \end{align}
    Adding \eqref{e:coco1} and \eqref{e:Lipds1} to \eqref{e:fire1a}, we obtain
\begin{align}\label{e:f01}
  &  \beta\|Cx_{k+1}-Cx\|^2\leq \frac{1}{\tau}\scal{x_k-x_{k+1}}{x_{k+1}-x}+\scal{(C+Q)x_{k+1}-(C+Q)x_k}{x_{k+1}-x}\notag\\
   & +\scal{\LL^*(\vv-\vv_k)}{x_{k+1}-x}+ \frac{\beta_0}{2}\big(\|x_k-x_{k-1}\|^2+\|x_k-x_{k+1}\|^2 \big)-\scal{Qx_k-Qx_{k-1}}{x_k-x},
\end{align}
which is equivalent to 
\begin{align}
  &  \beta\|Cx_{k+1}-Cx\|^2+\frac{1}{2\tau}\|x_{k+1}-x\|^2 \leq \frac{1}{2\tau}\|x_k-x\|^2+\scal{\LL^*(\vv-\vv_k)}{x_{k+1}-x}\notag\\
    &-\scal{Qx_k-Qx_{k-1}}{x_k-x}+\scal{(C+Q)x_{k+1}-(C+Q)x_k}{x_{k+1}-x}
    -\frac{1}{2\tau}\|x_k-x_{k+1}\|^2\notag\\
    &+\frac{\beta_0}{2}\big(\|x_k-x_{k-1}\|^2+\|x_k-x_{k+1}\|^2 \big). 
    \label{e:f02aa}
\end{align}
By denoting $t_k=\scal{Qx_k-Qx_{k-1}}{x_k-x}$, we can rewrite \eqref{e:f02aa} as 
\begin{align}
  &  \beta\|Cx_{k+1}-Cx\|^2+\frac{1}{2\tau}\|x_{k+1}-x\|^2 -t_{k+1}\leq \frac{1}{2\tau}\|x_k-x\|^2+\scal{\LL^*(\vv-\vv_k)}{x_{k+1}-x}\notag\\
    &-t_k+\scal{Cx_{k+1}-Cx_k}{x_{k+1}-x}
    -\frac{1}{2\tau}\|x_k-x_{k+1}\|^2
    +\frac{\beta_0}{2}\big(\|x_k-x_{k-1}\|^2+\|x_k-x_{k+1}\|^2 \big). 
    \label{e:f02}
\end{align}
For each $i\in\{1,\ldots,m\}$, define $w_{i,k}=v_{i,k}-\sigma_i L_{i}x_k$. Then, we derive from \eqref{algo:primal-dual} that
\begin{equation}
  \frac{1}{\sigma_i} (w_{i,k} - w_{i,k+1})+L_ix_{k+1} -D_{i}^{-1}v_{i,k} \in r_i+ B_{i}^{-1}(w_{i,k+1} +\sigma_i L_ix_k).\label{e:du1}
\end{equation}
Since $B_{i}^{-1}$ is monotone and $L_ix-D_{i}^{-1}v_i\in r_i+ B_{i}^{-1}v_i$, \eqref{e:split1} and \eqref{e:du1} yield 
\begin{align}
    0&\leq \scal{  \frac{1}{\sigma_i} (w_{i,k} - w_{i,k+1})+L_ix_{k+1}-L_ix-(D_{i}^{-1}v_{i,k}-D_{i}^{-1}v_i)}{w_{i,k+1} +\sigma_i L_ix_k-v_i}\notag\\
    &=\frac{1}{\sigma_i}\scal{w_{i,k} - w_{i,k+1}}{w_{i,k+1}+\sigma_i L_ix_k-v_i}
    +\scal{L_ix_{k+1}-L_ix}{w_{i,k+1} +\sigma_i L_ix_k-v_i}-\alpha_{i,k},\label{e:mid}
    \end{align}
    where we set 
    \begin{equation}
        \alpha_{i,k}= \scal{(D_{i}^{-1}v_{i,k}-D_{i}^{-1}v_i)}{w_{i,k+1} +\sigma_i L_ix_k-v_i}.
    \end{equation}
 Since $w_i= v_i-\sigma_i L_ix$, we can express the first term in \eqref{e:mid} as 
 \begin{align}
   & \frac{1}{\sigma_i}\scal{w_{i,k} - w_{i,k+1}}{w_{i,k+1}+\sigma_i L_ix_k-v_i} =    \frac{1}{\sigma_i}\scal{w_{i,k} - w_{i,k+1}}{w_{i,k+1}-w_i+\sigma_i L_i(x_k-x)}\notag\\
    &=  \frac{1}{\sigma_i}\scal{w_{i,k} - w_{i,k+1}}{w_{i,k+1}-w_i}
    +  \scal{v_{i,k} - v_{i,k+1}}{ L_i(x_k-x)}
    -\sigma_i\scal{L_i(x_i-x_{i+1})}{L_i(x_k-x)}.\notag
 \end{align}
 We can write the second term in \eqref{e:mid} as 
 \begin{align}
     \scal{L_ix_{k+1}-L_ix}{w_{i,k+1} +\sigma_i L_ix_k-v_i}&=
     \scal{L_ix_{k+1}-L_ix}{v_{i,k+1} -v_i}\notag\\
     &\hspace{1cm}+\sigma_i\scal{L_ix_{k+1}-L_ix}{ L_ix_k-L_ix_{k+1}}.
 \end{align}
 Therefore, \eqref{e:mid} becomes
    \begin{align}
           % 0&\leq \frac{1}{\sigma}\scal{w_k - w_{k+1}}{w_{k+1}-v}
   % + \scal{w_k - w_{k+1}}{Lx_k}
    % +\scal{Lx_{k+1}-Lx}{v_{k+1}-v}\notag\\
     %&\hfill+\sigma\scal{Lx_{k+1}-Lx}{Lx_k-Lx_{k+1}}-\alpha_k\notag\\
      % &=\frac{1}{\sigma}\scal{w_k - w_{k+1}}{w_{k+1}-w}
    %+ \scal{w_k - w_{k+1}}{Lx_k-Lx}
     %+\scal{Lx_{k+1}-Lx}{v_{k+1}-v}\notag\\
     %&\hfill+\sigma\scal{Lx_{k+1}-Lx}{Lx_k-Lx_{k+1}}-\alpha_k\notag\\
     0&\leq \frac{1}{\sigma_i}\scal{w_{i,k} - w_{i,k+1}}{w_{i,k+1}-w_i} +\scal{v_{i,k}-v_{i,k+1}}{L_i(x_k-x)} +\scal{L_i(x_{k+1}-x)}{v_{i,k+1}-v_i} \notag\\
     &-\sigma_i \|L_ix_k-L_ix_{k+1}\|^2-\alpha_{i,k}.
    \label{e:ff2}
\end{align}
We can express the first term in the right hand side of
\eqref{e:ff2} as follows
\begin{align}
  \frac{1}{\sigma_i}\scal{w_{i,k} - w_{i,k+1}}{w_{i,k+1}-w_i} 
  &=   \frac{1}{2\sigma_i}\|w_{i,k}-w_i\|^2
  -\frac{1}{2\sigma_i}\|w_{i,k+1}-w_i\|^2
  -\frac{1}{2\sigma_i}\|w_{i,k}-w_{i,k+1}\|^2\notag\\
    &=   \frac{1}{2\sigma_i}\|w_{i,k}-w_i\|^2
  -\frac{1}{2\sigma_i}\|w_{i,k+1}-w_i\|^2
  -\frac{1}{2\sigma_i}\|v_{i,k}-v_{i,k+1}\|^2\notag\\
  &\quad +\scal{v_{i,k}-v_{i,k+1}}{L_i(x_k-x_{k+1})} -\frac{\sigma_i}{2}\|L_ix_k-L_ix_{k+1}\|^2.\label{e:ff3}
\end{align}
Adding \eqref{e:ff2} and \eqref{e:ff3}, we obtain
\begin{align}
   \frac{1}{2\sigma_i}\|w_{i,k+1}-w_i\|^2 &\leq    \frac{1}{2\sigma_i}\|w_{i,k}-w_i\|^2 +\scal{v_{i,k}-v_{i,k+1}}{L_ix_k-L_ix} +\scal{L_ix_{k+1}-L_ix}{v_{i,k+1}-v_i} \notag\\
     &-\frac{3\sigma_i}{2} \|L_ix_k-L_ix_{k+1}\|^2  -\frac{1}{2\sigma_i}\|v_{i,k}-v_{i,k+1}\|^2+\scal{v_{i,k}-v_{i,k+1}}{L_ix_k-L_ix_{k+1}}\notag\\ 
     &-\alpha_{i,k}.
     \label{e:ff4}
\end{align}
Let us recall from \eqref{e:weigh} that
\begin{equation}
\begin{cases}
    \||\ww_k-\ww\||^2&=\frac{1}{2\tau}\|x_k-x\|^2+ \sum_{i=1}^m \frac{1}{2\sigma_i}\|w_{i,k}-w_i\|^2\; \text{where $\ww=(x,w)$, $\ww_k= (x_k,(w_{i,k})_{1\leq i\leq m})$}, \notag\\
   \||\xx_k-\xx_{k+1}\||^2&= \sum_{i=1}^m\big(  \frac{1}{2\sigma_i}\|v_{i,k}-v_{i,k+1}\|^2\big)+\frac{1}{2\tau}\|x_k-x_{k+1}\|^2\notag\\
    \end{cases}
\end{equation}
Then, by summing \eqref{e:ff4} over $i=1$ to $i=m$
and adding the result to \eqref{e:f02}, we obtain
\begin{align}
  \beta \|Cx_{k+1}-Cx\|^2 + & \||\ww_{k+1}-\ww\||^2 \leq    \||\ww_k-\ww\||^2-\||\xx_k-\xx_{k+1}\||^2
  +2\scal{\vv_k-\vv_{k+1}}{\LL(x_k-x_{k+1})}_{\GG}\notag\\
  &- \sum_{i=1}^m \frac{3\sigma_i}{2}\|L_i(x_k-x_{k+1})\|^2 -\alpha_k+\scal{Cx_{k+1}-Cx_k}{x_{k+1}-x} \notag\\
 & +\frac{\beta_0}{2}\big(\|x_k-x_{k-1}\|^2+\|x_k-x_{k+1}\|^2\big) +t_{k+1}-t_k,\label{e:ff0a}
\end{align}
where we set 
\begin{equation}
    \alpha_k= \sum_{i=1}^m \alpha_{i,k}.
\end{equation}
We next express the last term in \eqref{e:ff0a} as
\begin{align}
    \scal{Cx_{k+1}-Cx_k}{x_{k+1}-x} &= 
    \scal{Cx_{k+1}-Cx}{x_{k+1}-x}+\scal{Cx-Cx_k}{x_{k+1}-x}\notag\\
    &=\scal{Cx_{k+1}-Cx}{x_{k+1}-x}+\scal{Cx-Cx_k}{x_{k}-x}\notag\\
    &+ \scal{Cx-Cx_k}{ x_{k+1}-x_k}\notag\\
    &\leq \scal{Cx_{k+1}-Cx}{x_{k+1}-x}+\scal{Cx-Cx_k}{x_{k}-x}\notag\\
    &\quad+\beta' \|Cx-Cx_k\|^2+ \frac{1}{4\beta'}\| x_{k+1}-x_k\|^2.\label{e:CoC}
\end{align}
Since $D_{i}^{-1}$ is $\mu_i$-cocoercive, we can express $\alpha_{i,k}$ as follows: 
\begin{align}\label{e:sumal}
   - \alpha_{i,k}&= -\scal{D_{i}^{-1}v_{i,k}-D_{i}^{-1}v_i}{w_{i,k+1}+\sigma_i L_ix_k-v_i}\notag\\
    &= \scal{D_{i}^{-1}v-D_{i}^{-1}v_{i,k}}{v_{i,k+1}-v_i} +
    \sigma_i \scal{D_{i}^{-1}v_i-D_{i}^{-1}v_{i,k}}{L_ix_k-L_ix_{k+1}}\notag\\
    &= \scal{D_{i}^{-1}v_i-D_{i}^{-1}v_{i,k+1}}{v_{i,k+1}-v_i} +\scal{D_{i}^{-1}v_{i,k+1}-D_{i}^{-1}v_{i,k}}{v_{i,k+1}-v_i}\notag\\
    &\hspace{5.5cm}+\sigma_i \scal{D_{i}^{-1}v_i-D_{i}^{-1}v_{i,k}}{L_ix_k-L_ix_{k+1}}\notag\\
    &\leq -\mu_i \|D_{i}^{-1}v_{i,k+1}-D_{i}^{-1}v_i\|^2 
    +\scal{D_{i}^{-1}v_{i,k+1}-D_{i}^{-1}v_i}{v_{i,k+1}-v_i} \notag\\
    &+\scal{D_{i}^{-1}v_i-D_{i}^{-1}v_{i,k}}{v_{i,k}-v_i}
    +\sigma_i \scal{D_{i}^{-1}v_i-D^{-1}_{i}v_{i,k}}{L_ix_k-L_ix_{k+1}}\notag\\
   & \hspace{6cm}+\scal{D_{i}^{-1}v_i-D^{-1}_iv_{i,k}}{v_{i,k+1}-v_{i,k}}\notag\\
    &\leq 
    \scal{D^{-1}_{i}v_{i,k+1}-D^{-1}_iv_i}{v_{i,k+1}-v_i} +\scal{D^{-1}_iv_i-D^{-1}_iv_{i,k}}{v_{i,k}-v_i}\notag\\
    &+\mu'_i \|D^{-1}_iv_{i,k}-D^{-1}_iv_i\|^2 -\mu_i \|D_{i}^{-1}v_{i,k+1}-D_{i}^{-1}v_i\|^2 +\frac{\sigma_{i}^2}{2\mu'_i}\|L_ix_k-L_ix_{k+1}\|^2\notag\\
    &\hspace{8cm}+\frac{1}{2\mu'_i}\|v_{i,k}-v_{i,k+1}\|^2.
\end{align}
Let us set 
\begin{equation}
    \omega_k= \beta' \|Cx_k-Cx\|^2+\sum_{i=1}^m\mu'_i\|D_{i}^{-1}v_{i,k}-D_{i}^{-1}v_i\|^2.
\end{equation}
Now, summing \eqref{e:sumal} over $i=1$ to $i=m$ and adding the result to \eqref{e:CoC}, we get
\begin{align}
  -\alpha_k&+\scal{Cx_{k+1}-Cx_k}{x_{k+1}-x}= 
  \scal{\QQ\xx_{k+1}-\QQ\xx}{\xx_{k+1}-\xx}_{\KK}
  +\scal{\QQ\xx-\QQ\xx_k }{\xx_k-\xx}_{\KK}\notag\\
  &+\omega_k-\omega_{k+1}+\frac{1}{4\beta'}\|x_{k+1}-x_k\|^2
  +\sum_{i=1}^m\bigg( \frac{\sigma_{i}^2}{2\mu'_i}\|L_ix_k-L_ix_{k+1}\|^2
    +\frac{1}{2\mu'_i}\|v_{i,k}-v_{i,k+1}\|^2 \bigg)\notag\\
    &+\beta' \|Cx_{k+1}-Cx\|^2
    +\sum_{i=1}^m(\mu'_i-\mu_i)\|D_{i}^{-1}v_{i,k+1}-D_{i}^{-1}v_i\|^2.\notag
\end{align}
Therefore, we further estimate \eqref{e:ff0a} as
\begin{align}
 &\||\ww_{k+1}-\ww\||^2 +\omega_{k+1} +
 \scal{\QQ\xx-\QQ\xx_{k+1}}{\xx_{k+1}-\xx}_{\KK}
 \leq    \||\ww_k-\ww\||^2 + \omega_k
  +\scal{\QQ\xx-\QQ\xx_k }{\xx_k-\xx}_{\KK}\notag\\
 &-\||\xx_k-\xx_{k+1}\||^2
  - \sum_{i=1}^m\sigma_i( \frac{3}{2}-\frac{\sigma_i}{2\mu'_i})\|L_i(x_k-x_{k+1})\|^2+\frac{1}{4\beta'} \|x_{k+1}-x_k\|^2
  +\sum_{i=1}^m\frac{1}{2\mu'_i}\|v_{i,k}-v_{i,k+1}\|^2 \notag\\ &+2\scal{\vv_k-\vv_{k+1}}{\LL(x_k-x_{k+1})}_{\GG}-e_k+\frac{\beta_0}{2}\big(\|x_k-x_{k-1}\|^2+\|x_k-x_{k+1}\|^2\big) +t_{k+1}-t_k,
  \label{e:ff0a1}
\end{align}
where we set 
\begin{equation}
    e_k = (\beta-\beta') \|Cx_{k+1}-Cx\|^2
    +\sum_{i=1}^m(\mu_i-\mu'_i)\|D_{i}^{-1}v_{i,k+1}-D_{i}^{-1}v_i\|^2.
\end{equation}
Using the Cauchy Schwarz inequality for the last term in \eqref{e:ff0a1}, we have
\begin{align}
  &2\scal{\vv_k-\vv_{k+1}}{\LL x_k-\LL x_{k+1}} 
  \leq \sum_{i=1}^m 2\|v_{i,k}-v_{i,k+1}\| \|L_ix_k-L_ix_{k+1}\|  \notag\\
 &\leq \sum_{i=1}^m (\frac{1}{2\sigma_i}-\frac{1}{2\mu'_i} -\frac{\epsilon}{2\sigma_i\mu'_i})\|v_{i,k}-v_{i,k+1}\|^2+ \sum_{i=1}^m(\frac{1}{2\sigma_i}-\frac{1}{2\mu'_i} -\frac{\epsilon}{2\sigma_i\mu'_i})^{-1} \|L_ix_k-L_ix_{k+1}\|^2,
\end{align}
where $\epsilon \in \left]0,\min_{1\leq i\leq m}(\mu'_i-\sigma_i)\right[$. Hence, we derive from \eqref{e:ff0a1} that
\begin{align}
 &\||\ww_{k+1}-\ww\||^2 +\omega_{k+1} +
 \scal{\QQ\xx-\QQ\xx_{k+1}}{\xx_{k+1}-\xx}_{\KK}
 \leq    \||\ww_k-\ww\||^2 + \omega_k
  +\scal{\QQ\xx-\QQ\xx_k }{\xx_k-\xx}_{\KK}\notag\\
 &-\Big(\frac{1}{2\tau}-\frac{1}{4\beta'}\Big)\|x_k-x_{k+1}\|^2
  + \sum_{i=1}^m\sigma_i\Big( \frac{\sigma_i}{2\mu'_i}+\frac{2\mu'_i}{\mu'_i-\sigma_i-\epsilon}-\frac{3}{2}\Big)\|L_i(x_k-x_{k+1})\|^2\notag\\
  &-\sum_{i=1}^m\frac{\epsilon}{2\sigma_i\mu'_i}\|v_{i,k}-v_{i,k+1}\|^2-e_k.
  \label{e:ff0a31}
\end{align}
Observe that 
\begin{equation}
  \Big(\forall \epsilon \in ]0,\min_{1\leq i\leq m}(\mu'_i-\sigma_i)[\Big)\;  \frac{\sigma_i}{2\mu'_i}+\frac{2\mu'_i}{\mu'_i-\sigma_i-\epsilon}-\frac{3}{2} > 0, 
\end{equation}
and hence
\begin{align}
 &\||\ww_{k+1}-\ww\||^2 +\omega_{k+1} +
 \scal{\QQ\xx-\QQ\xx_{k+1}}{\xx_{k+1}-\xx}_{\KK}-t_{k+1} + \frac{\beta_0}{2}\|x_k-x_{k+1}\|^2
 \leq  \notag\\
& \||\ww_k-\ww\||^2 + \omega_k
  +\scal{\QQ\xx-\QQ\xx_k }{\xx_k-\xx}_{\KK}-t_k+ \frac{\beta_0}{2}\|x_k-x_{k-1}\|^2\notag\\
 &-\bigg(\frac{1}{2\tau}-\frac{1}{4\beta'}-\beta_0
  - \sum_{i=1}^m\sigma_i\Big( \frac{\sigma_i}{2\mu'_i}+\frac{2\mu'_i}{\mu'_i-\sigma_i-\epsilon}-\frac{3}{2}\Big)\|L_i\|^2\bigg) \|x_k-x_{k+1}\|^2\notag\\
  &-\sum_{i=1}^m\frac{\epsilon}{2\sigma_i\mu'_i}\|v_{i,k}-v_{i,k+1}\|^2-e_k.
  \label{e:ff0a32}
\end{align}
Moreover, it follows from \eqref{e:constep1} that
\begin{equation}
  \Big(\exists \epsilon \in ]0,\min_{1\leq i\leq m}(\mu'_i-\sigma_i)[\Big)\; \frac{1}{2\tau}-\frac{1}{4\beta'}- \beta_0-\sum_{i=1}^m\sigma_i\|L_i\|^2\Big(\frac{\sigma_i}{2\mu'_i}+\frac{2\mu'_i}{\mu'_i-\sigma_i-\epsilon}-\frac{3}{2}\Big) \geq \frac{\epsilon}{2\tau \max\mu'_i}.
\end{equation}
Thus, the estimation \eqref{e:ff0a32} is simplified to 
\begin{align}
&\||\ww_{k+1}-\ww\||^2 +\omega_{k+1} +
 \scal{\QQ\xx-\QQ\xx_{k+1}}{\xx_{k+1}-\xx}_{\KK}-t_{k+1} + \frac{\beta_0}{2}\|x_k-x_{k+1}\|^2
 \leq \notag\\
& \||\ww_k-\ww\||^2 + \omega_k
  +\scal{\QQ\xx-\QQ\xx_k }{\xx_k-\xx}_{\KK}-t_k+ \frac{\beta_0}{2}\|x_k-x_{k-1}\|^2\notag\\
- &\frac{\epsilon}{\max \mu'_i}\||\xx_{k+1}-\xx_k\||^2-e_k.
  \label{e:ff0a3}
\end{align}
Next, we find a lower bound of  
\begin{equation}
\Gamma_k:=\||\ww_k-\ww\||^2 + \omega_k
  +\scal{\QQ\xx-\QQ\xx_k }{\xx_k-\xx}_{\KK}-t_k+ \frac{\beta_0}{2}\|x_k-x_{k-1}\|^2.
\end{equation}
Using the Cauchy-Schwarz inequality, we have the following.
  \begin{align}
  \Big(\forall i\in\{1,\ldots,m\}\Big)\;
  \begin{cases}
     & \beta' \|Cx_k-Cx\|^2 +\scal{Cx-Cx_k}{x_k-x}\geq -\frac{1}{4\beta'}\|x_k-x\|^2\notag\\
     & \mu'_i\|D_{i}^{-1}v_{i,k}-D^{-1}_i v_i\|^2+\scal{D_{i}^{-1}v_{i,k}-D^{-1}_i v_i}{v_{i,k}-v_i}
     \geq -\frac{1}{4\mu'_i} \|v_{i,k}-v_i\|^2,
      \end{cases}
  \end{align}
  and
  \begin{align}
      -t_k+ \frac{\beta_0}{2}\|x_k-x_{k-1}\|^2&= -\scal{Qx_k-Qx_{k-1}}{x_k-x}+\frac{\beta_0}{2}\|x_k-x_{k-1}\|^2\notag\\
     & \geq -\|Qx_k-Qx_{k-1}\|\|x_k-x\|+\frac{\beta_0}{2}\|x_k-x_{k-1}\|^2\notag\\
      &\geq -\frac{\beta_0}{2}\|x_k-x\|^2.
  \end{align}
  It follows that 
  \begin{equation}
   \omega_k
  +\scal{\QQ\xx-\QQ\xx_k }{\xx_k-\xx}_{\KK}-t_k+ \frac{\beta_0}{2}\|x_k-x_{k-1}\|^2\geq 
  - (\frac{1}{4\beta'}+\frac{\beta_0}{2})\|x_k-x\|^2-\sum_{i=1}^m\frac{1}{4\mu'_i} \|v_{i,k}-v_i\|^2.\label{e:dsh}
  \end{equation}
  Adding $\||\ww_k-\ww\||^2$ to both sides of \eqref{e:dsh}, we obtain
  \begin{align}
  \Gamma_k
 \geq \Big(\frac{1}{2\tau} - \frac{1}{4\beta'}-\frac{\beta_0}{2}\Big)\|x_k-x\|^2+\sum_{i=1}^m\Big(\frac{1}{2\sigma_i} \|w_{i,k}-w_i\|^2- \frac{1}{4\mu'_i} \|v_{i,k}-v_i\|^2\Big).\label{e:dshs}  
  \end{align}
 For each $i\in\{1,\ldots,m\}$, by the definitions of $w_{i,k} = v_{i,k}-\sigma_i L_ix_k$ and $w_i =  v_i- \sigma_i L_ix$, we have
  \begin{align}
    \frac{1}{2\sigma_i}\|w_{i,k}-w_i\|^2 - \frac{1}{4\mu'_i} \|v_{i,k}-v_i\|^2 &= (\frac{1}{2\sigma_i} -\frac{1}{4\mu'_i})\|v_{i,k}-v_i\|^2 \notag\\
    &+\frac{\sigma_i}{2}\|L_ix_k-L_ix\|^2 -\scal{v_{i,k}-v_i}{L_ix_k-L_ix}\notag\\
    &\geq \frac{1}{4\sigma_i} \|v_{i,k}-v_i\|^2
    +(\frac{\sigma_i}{2} - \frac{\sigma_i\mu'_i}{\mu'_i-\sigma_i})\|L_ix_k-L_ix\|^2\notag\\
      &\geq \frac{1}{4\sigma_i} \|v_{i,k}-v_i\|^2
    +(\frac{\sigma_i}{2} - \frac{\sigma_i\mu'_i}{\mu'_i-\sigma_i})\|L_i\|^2\|x_k-x\|^2.
 \end{align}
 In turn, 
   \begin{align}
  \Gamma_k&\geq \Big(\frac{1}{2\tau} - \frac{1}{4\beta'} - \frac{\beta_0}{2}-\sum_{i=1}^m \Big( \frac{\sigma_i\mu'_i}{\mu'_i-\sigma_i}-\frac{\sigma_i}{2}\Big)\|L_i\|^2\Big)\|x_k-x\|^2+\sum_{i=1}^m\frac{1}{4\sigma_i} \|v_{i,k}-v_i\|^2\notag\\
 &\geq \frac{\epsilon}{2\tau\max\mu'_i}\|x_k-x\|^2
 +\sum_{i=1}^m\frac{1}{4\sigma_i} \|v_{i,k}-v_i\|^2\notag\\
 &\geq 0.
 \label{e:dshs}  
  \end{align}
\ref{t:main1i}: We derive from \eqref{e:ff0a3} and Lemma \ref{lem:polyak} that
\begin{equation}\label{e:pro1}
    \begin{cases}
        &\exists \lim \||\ww_k-\ww\||^2+   \omega_k
  +\scal{\QQ\xx-\QQ\xx_k }{\xx_k-\xx}_{\KK}-t_k+\frac{\beta_0}{2}\|x_k-x_{k-1}\|^2 =\overline{\xi} \in [0,+\infty[,\\
  &\sum_{k\in\NN}\||\xx_{k+1}-\xx_k\||^2 < +\infty \;
  \text{and}\; \sum_{k\in\NN} e_k < +\infty.
    \end{cases}   
\end{equation}
In particular, we obtain 
\begin{equation}\label{e:pro2}
    \sum_{k\in\NN}\|x_{k+1}-x_k\|^2 < +\infty
   \; \text{and}\; (\forall i\in\{1,\ldots,m\})\; \sum_{k\in\NN}\|v_{i,k+1}-v_{i,k}\|^2<+\infty,
\end{equation}
that proves \ref{t:main1i}. 
Let us prove \ref{t:main1ii} $\&$ \ref{t:main1iii}. It follows from \eqref{e:pro1} that 
\begin{equation}
    \lim \||\ww_k-\ww\||^2 = \overline{\xi},\label{e:conv1}
\end{equation}
hence $(\ww_k)_{k\in\NN}$ is bounded. Moreover, \eqref{e:dshs} implies that the sequence $(x_k,v_{1,k},\ldots, v_{m,k})_{k\in\NN}$ is also bounded. 
Let $\yy^*=(y^*,w_{1}^*,\ldots,w_{m}^*)$ be a weak
cluster point of $(\ww_k)_{k\in\NN}$, i.e., there exists a subsequence 
$(x_{n_k},w_{1,n_k},\ldots, w_{m,n_k})_{k\in\NN}$ that converges weakly to $\yy^*$. It follows that 
$(x_{n_k},v_{1,n_k},\ldots, v_{m,n_k})_{k\in\NN}$ converges weakly to 
$(y^*,w_{1}^*+\sigma_1 L_1 y^*, \ldots, w_{m}^*+\sigma_m L_m y^*).$ 
By \ref{t:main1i}, the sequence
$(x_{1+n_k},v_{1,1+n_k},\ldots, v_{m,1+n_k})_{k\in\NN}$  converges weakly to $(y^*,w_{1}^*+\sigma_1 L_1 y^*, \ldots, w_{m}^*+\sigma_m L_m y^*)$. We next use \eqref{e:priup1} and \eqref{e:du1},
\begin{equation}
    \begin{cases}
        &\frac{x_{n_k}-x_{n_k+1}}{\tau} - Cx_{n_k}+ Cx_{1+n_k}-2Qx_{n_k}+Qx_{n_k-1} +Qx_{n_k+1}-\LL^*\vv_{n_k}+\sum_{i=1}^mL_{i}^*(w_{i,1+n_k}+\sigma_i L_i x_{n_k}) \\
        &\hspace{4cm} \in -z+(A+C+Q)x_{n_k+1}+\sum_{i=1}^mL_{i}^*(w_{i,1+n_k}+\sigma_i L_i x_{n_k}) \\
        &\frac{1}{\sigma_i} (w_{i,n_k} - w_{i,n_k+1})-D_{i}^{-1}v_{i,n_k} + D_{i}^{-1}(w_{i,n_k+1} +\sigma_i L_ix_{n_k}) \\
       &\hspace{4cm} \in r_i+(B_{i}^{-1}+D_{i}^{-1})(w_{i,n_k+1} +\sigma_i L_ix_{n_k})-L_ix_{n_k+1}.\label{e:cbc1}
    \end{cases}
    \end{equation}
    Note that the sequence on the left hand side of \eqref{e:cbc1} converges strongly to $0$.
Set $(\forall k\in\NN)\; \af_k = (x_{k+1},\overline{\ww}_{k})$ where $\overline{\ww}_k= \ww_{1+k}+(\sigma_i L_ix_{k})_{1\leq i\leq m}$. Then, $\af_{n_k}\weakly \af^*=(y^*,w_{1}^*+\sigma_1 L_1 y^*, \ldots, w_{m}^*+\sigma_m L_m y^*)$.
Then the right hand side of \eqref{e:cbc1} is $(\MM+\QQ+\SSS)\af_{n_k}$. Since $A, C$ and $(B_i)_{1\leq i\leq m}$, $(D_i)_{1\leq i\leq m}$ are maximally monotone, $\MM$ and $\QQ$ are also maximally monotone \cite[Proposition 20.22, Proposition 23.18 ]{livre1}. Since $\SSS$ is maximally monotone and $\SSS+\QQ$ is Lipschitz, the sum $\MM+\QQ+\SSS$ is maximally monotone \cite[Corollary 25.5]{livre1}.  
Hence, its graph is closed in $\HH^{weak}\times\HH^{strong}$ \cite[Proposition 20.38]{livre1}. It follows that $\af^*\in \zer(\MM+\QQ+\SSS)$ and hence $\yy^*\in\overline{Z}$. We conclude that every weak cluster point of 
$(\ww_k)_{k\in\NN}$ is in $\overline{Z}$. Combining this with \eqref{e:conv1}, by Lemma  \cite[Lemma 2.47]{livre1}, the whole sequence $(\ww_k)_{k\in\NN}$ converges weakly to some point $\ww=(\overline{x},\overline{v}_1-\sigma_1 L_1\overline{x},\ldots, \overline{v}_m-\sigma_m L_m\overline{x})$ such that
$(\overline{x},\overline{v}_1,\ldots,\overline{v}_m)\in Z$. Moreover,
$(x_k,v_{1,k},\ldots, v_{m,k})\weakly (\overline{x},\overline{v}_1,\ldots,\overline{v}_m)$.
Lemma \ref{l:1} yields $\overline{x}\in\mathcal{P} $ and 
$(\overline{v}_1,\ldots,\overline{v}_m)\in\mathcal{D}$.\\
\ref{t:main1ii}: The result follows from \eqref{e:pro2} and $\overline{x}\in\mathcal{P}$, and 
$(\overline{v}_1,\ldots,\overline{v}_m)\in\mathcal{D}$.\\
\ref{t:mainiv}: Suppose that $C$ is uniformly monotone at $\overline{x}$. Then there exists an increasing function $\phi_C\colon \left[0,+\infty\right[\to\left[0,+\infty\right]$ vanishing only at $0$ such that 
\begin{equation}
    \phi_C(\|x_k-\overline{x}\|)\leq \scal{x_k-\overline{x}}{Cx_k-C\overline{x}}
    \leq \|x_k-\overline{x}\|\|Cx_k -C\overline{x}\|.
\end{equation}
By \ref{t:maini2}, $Cx_k\to C\overline{x}$ and by \ref{t:main1ii}, $(\|x_k-\overline{x}\|)_{k\in\NN}$ is bounded. It follows that $\phi_C(\|x_k-\overline{x}\|)\to 0$ and hence $x_k\to \overline{x}$. Now suppose that $A$ is uniformly monotone at $x=\overline{x}\in \mathcal{P}$. Then, \eqref{e:fire1a} is improved to 
  \begin{align}
   \phi_A(\|x_{k+1}-x\|)&\leq
     \frac{1}{\tau}\scal{x_k-x_{k+1}}{x_{k+1}-x}+\scal{Cx-Cx_k}{x_{k+1}-x} +\scal{\LL^*\vv-\LL^*\vv_k}{x_{k+1}-x}\notag\\
   &\quad +\scal{Qx-Qx_k}{x_{k+1}-x}-\scal{Qx_k-Qx_{k-1}}{x_{k+1}-x}.\notag
\end{align}
Hence, \eqref{e:ff0a3} is improved to
\begin{align}
 \Gamma_{k+1} \leq \Gamma_k
- \frac{\epsilon}{\max \mu'_i}\||\xx_{k+1}-\xx_k\||^2-e_k- \phi_A(\|x_{k+1}-x\|),
  \label{e:ff0a3aa}
\end{align}
which implies that $\sum_{k\in\NN} \phi_A(\|x_{k+1}-x\|)<+\infty$. Therefore, $x_k\to x=\overline{x}.$\\
\ref{t:mainv}: Using the same fashion as the proof of \ref{t:mainiv}.
\end{proof}
\begin{remark} Here are some comments.
\begin{itemize}
\item[(i)]The main feature of the proposed algorithm is that the term $\sigma_i(L_ix_{k+1}-L_ix_k)$ is not involved in the resolvent $J_{\sigma_iB_{i}^{-1}}$. This idea  first appeared in \cite{Tam19} where the shadow Douglas-Rachford was introduced. Hence, we refer \eqref{algo:primal-dual} as a primal dual shadow Douglas-Rachford splitting algorithm.
\item[(ii)]
The algorithm \eqref{algo:primal-dual} differs from all existing primal--dual methods for Problem \ref{p:Prob0} available in the literature; see, for example, \cite{acta2,Buicomb,Condat,vu2013,Icip14,plc6,Bot2013,Bot2014}. 
\item[(iii)] When $m=1$ and $D_{1}^{-1}=0$, alternative methods can be found in \cite{Pontus23}. In addition, when $C=0$, the condition \eqref{e:constep} becomes
\begin{equation}
    \frac{1}{2\tau}
 -\beta_0 - \frac{1}{2}\sigma_1\|L_1\|^2>0,
\end{equation}
which is the same as the condition ensuring the weak convergence of the alternative method proposed in
 \cite{Malistky26}.
\item[(iv)]  When $m=1$ and $D_{1}^{-1}=0, Q=0$, the algorithm \eqref{algo:primal-dual} is equivalent to the method in \cite[Eq. (38)]{Dao26} which itself is an
extension of the primal-dual algorithm in \cite{Tam19} which is designed to solve the  bilinear convex-concave saddle point problem.  In this case, Problem \ref{p:Prob0} can be solved by the split-Douglas-Rachford  in \cite{Luis21}.
 \item[(iv)] When $m=1$ and $D_{1}^{-1}=0$, the condition on stepsize \eqref{e:constep} is exactly the same as that of algorithm in \cite{Condat}, that is 
 \begin{equation}
     \frac{1}{\tau} -\sigma_1\|L_1\|^2 > \frac{1}{2\beta}.
 \end{equation}
However, in the general case where $m>1$ and $\exists D_{i}^{-1}\neq 0$, the condition \eqref{e:constep} differs from that in \cite{vu2013}, while yielding the same type of convergence results.
\end{itemize}
\end{remark}
\section{Application to structured minimization problems}
\label{s:O-inf}
In this section, we focus on an application of the proposed algorithm to solve structured minimization problems. Let us recall 
several notions in convex optimization \cite{livre1}. We denote by $\Gamma_0(\HH)$ the class of lower semicontinuous convex functions $f\colon \HH\to\left]-\infty,+\infty\right]$ such that
$
\dom(f)=\{x\in \HH\mid f(x)<+\infty\}\neq\varnothing.
$
Let $f\in\Gamma_0(\HH)$. The subdifferential of $f$ at $x \in \dom f$ is $$ \partial f(x)=\{u \in \HH,\ \forall z \in \dom f:\ f(z) \ge f(x)+\scal{u}{z-x} \}.$$
The proximity operator of $f$ is
\begin{align}
\prox_f:\mathcal{H} \to \mathcal{H}: x \mapsto \underset{y \in \mathcal{H}}{\operatorname{argmin}}\big(f(y)+\frac{1}{2}\|x-y\|^2\big).
\end{align}
 The conjugate function of $f$ is
\begin{align}
f^*:a \mapsto \underset{x \in \mathcal{H}}{\sup} \big( \scal{a}{x}-f(x) \big).
\end{align}
Note that $ (\partial f)^{-1}=\partial f^*$, and $p=\prox_{f}x$ iff $x-p\in\partial f(p)$. 
%The infimal convolution of the two functions $\ell$ and $g$ from $\HH$ to $\left]-\infty,+\infty\right]$ is 
 %\begin{equation}
 %\ell\;\vuo\; g\colon x \mapsto \inf_{y\in\HH}(\ell(y)+g(x-y)).
%\end{equation}
%Let $L\colon\HH\to\GG$ be a bounded linear operator. The infimal postcomposition of $f$ by $L$ is 
%\begin{equation}
%L\triangleright f:\GG\to[-\infty,+\infty]\colon
%y\mapsto \inf_{\substack{x\in\HH\\\substack{Lx=y}}} f(x).
%\end{equation}
Finally, let $C$ be a convex subset of $\HH$, the strong relative interior of $C$ is denoted by $\sri C$.
\begin{lemma}\label{L:repha} Under the setting of Problem \ref{prob:app0}. We have, 
\begin{equation}\label{e:vuoO}
    (\forall i\in \{1,\ldots,m\})\;  g_i\;\vuoO_i\;(W_i \triangleright\ell_i)= (O_i \triangleright g_i)\;\vuo\;(W_i \triangleright\ell_i).
\end{equation}
  Consequently, 
  \begin{equation}\label{e:Conj}
    (\forall i\in \{1,\ldots,m\})\; ( g_i\;\vuoO_i\;(W_i \triangleright\ell_i))^*= ((O_i \triangleright g_i)\;\vuo\;(W_i \triangleright\ell_i))^* = (O_i \triangleright g_i)^* +(W_i \triangleright \ell_i)^*.
\end{equation}
\end{lemma}
\begin{proof}Let $i\in \{1,\ldots,m\}$ and set $\ell =W_i \triangleright\ell_i$, and $w\in \GG_i$. We have
\begin{align}
       \big( (O_i \triangleright g_i)\;\vuo\;(W_i \triangleright\ell_i)\big)(w) &= \inf_{y\in \GG_i} [(O_i \triangleright g_i)(y)+\ell(w-y)]\notag\\
       &= \inf_{y\in \GG_i}[\inf_{\substack{u \in \YY_i \\ O_iu = y}} g_i(u) +\ell(w-y)]\notag\\
       &= \inf_{\substack{u \in \YY_i }} g_i(u) +\ell(w-O_iu)\notag\\
       &= (g_i\;\vuoO_i \ell)(w),
\end{align}
which proves \eqref{e:vuoO} and hence \eqref{e:Conj} follows directly.
\end{proof}

\noindent Lemma \ref{L:repha} allows us to rewrite the $O_i$-infimal convolution in terms of the infimal postcomposition. Therefore, we can now restate the problem of interest as follows.
\begin{problem}
\label{prob:app}
Under the setting of Problem \ref{prob:app0}. The primal problem \eqref{e:pri10} can be restated as
\begin{equation}\label{e:pri1}
\underset{x\in\HH}{\text{minimize}} \;
f(x)+\sum_{i=1}^m\big((O_i \triangleright g_i)\;\vuo\;(W_i \triangleright\ell_i)\big)(L_ix-r_i) + h(x) -\scal{z}{x},
\end{equation}
and its dual problem  is to
\begin{equation}\label{e:dual1}
 \underset{v_1\in\GG_1,\ldots,v_m\in\GG_m}{\text{minimize}} \;
(f^*\;\vuo\; h^*)\bigg(z-\sum_{i=1}^mL_{i}^*v_i\bigg) +
\sum_{i=1}^m \big((O_i \triangleright g_i)^*(v_i)+(W_i \triangleright\ell_i)^*(v_i) +\scal{v_i}{r_i}\big).
\end{equation}
We denote by $\mathcal{P}_1$ and $\mathcal{D}_1$ the sets of 
solutions to~\eqref{e:pri1} and~\eqref{e:dual1}, respectively.
\end{problem}
We note that, in the special case where
$$
(\forall i\in\{1,\ldots,m\})\quad
\XX_i=\YY_i=\GG_i
\quad\text{and}\quad
O_i=W_i=\Id,
$$
Problem \ref{prob:app} reduces to the problem considered in \cite{plc6}, which itself encompasses a broad class of problems arising in image processing \cite{cp16,pq15}, statistical learning \cite{livre2}, machine learning \cite{Bian24}, transport optimization \cite{Fadili}. 
\begin{lemma}\label{l:lips} Let $\ell \in \Gamma_0(\XX)$ be  $\mu$-strong convex 
for some $\mu \in \left]0,+\infty\right[$,
and $W\colon \XX\to\GG$ be a bounded linear operator. Then $(W\triangleright \ell)^*$ is differentiable with $\mu^{-1}\|W\|^2$-Lipschitz gradient, $\nabla (W\triangleright \ell)^*=W\circ\nabla\ell^*\circ W^*$. 
\end{lemma}
\begin{proof} Since $\ell$ is $\mu$-strongly convex, its conjugate $\ell^*$ is is Fr\'echet differentiable with $1/\mu$-Lipschitz gradient. Hence, $\dom\ell^* =\XX$ which implies that the condition $0\in\sri(\ran W^*-\dom \ell)$ is trivially satisified. In turn, $W\triangleright \ell \in \Gamma_0(\GG)$ and by \cite[Proposition 13.24]{livre1},  $(W\triangleright \ell)^*=\ell^*\circ W^*$ which is differentiable with $\mu^{-1}\|W\|^2$ gradient.
\end{proof}
\begin{remark}\label{Re:1} Let us set $\GG= \bigoplus_{i=1}^m\GG_i$. Let $\vv=(v_1,\ldots,v_m)$ be the generic element of $\GG$. Define, for each $i\in\{1,\ldots,m\}$, the projection operator $\pi_i\colon\GG\to\GG_i\colon \vv\mapsto v_i$. Define
\begin{equation}
\begin{cases}
    h_0&\colon \vv \mapsto \sum_{i=1}^m [(W_i \triangleright\ell_i)^*(v_i)+\scal{v_i}{r_i}]\\
    K_0&\colon \vv\mapsto -\sum_{i=1}^mL_{i}^*v_i\\
    K_i&\colon \vv\mapsto O_{i}^*\circ\pi_i\vv.
\end{cases}
\end{equation} We set
\begin{equation}
d_{0}  =f^*, l_0 =h^*, d_{i} = g_{i}^*, 
l_i=\iota_{\{0\}}, r^{*}_0=-z, r^{*}_i=0,
\text{where $i\in\{1,\ldots,m\}$}. 
\end{equation}
Then, the dual problem \eqref{e:dual1} can be reformulated as follows 
    \begin{equation}\label{e:Rdual}
 \underset{\vv\in\GG}{\text{minimize}} \; h_0(\vv) +\sum_{i=0}^m (d_i\vuo l_i)(K_i\vv-r^{*}_i).
\end{equation}
\end{remark}
\begin{theorem} \label{t:main2}
In Problem \ref{prob:app}, suppose that 
\begin{equation}\label{e:ran}
    z\in \ran\Bigg(\partial f + \sum_{i=1}^mL_{i}^*\big(\partial(O_i \triangleright g_i)\;\vuo\;\partial(W_i \triangleright\ell_i)\big)(L_i\cdot-r_i) + \nabla h\Bigg).
\end{equation}
Let $\tau$ and $(\sigma_i)_{0\leq i\leq m}$ be strictly positive real numbers such that 
\begin{equation}\label{e:constep1}
    \begin{cases}
   \sigma_0 < \beta,\\
        \frac{1}{2\tau}-\max_{1\leq i\leq m}\mu_{i}^{-1}\|W_i\|^2
 - \sum_{i=1}^m\|L_i\|^2\sigma_0\Big( \frac{\sigma_0}{2\beta}+\frac{2\beta}{\beta-\sigma_0}-\frac{3}{2}\Big) 
 -\sum_{i=1}^m \frac{\sigma_i}{2}\|O_i\|^2
 > 0.
    \end{cases}
\end{equation}
Let $x_{0}\in \HH$. For each $i\in\{1,\ldots,m\}$, let $w_{i,0}\in \YY_i$ and $v_{i,0}\in\GG_i$. Iterate
\begin{equation}
\label{algo:apprimal-dual0}
\begin{array}{l}
\left\lfloor
\begin{array}{l}
\operatorname{For \; \text{$i =1,\ldots,m$}}\\
\begin{array}{l}
\left\lfloor
\begin{array}{l}
v_{i,k+1} = v_{i,k}-\tau \big(O_iw_{i,k}-L_ix_k+ 
W_i\nabla\ell^{*}_i(W_{i}^*v_{i,k})-r_i\big)\\
\end{array}
\right.\\[2mm]
\end{array}\\
x_{k+1} = \prox_{\sigma_0 f}\Big(x_k-\sigma_0 \big(\sum_{i=1}^m L_{i}^*v_{i,k+1}-\nabla h(x_k)+z\big) \Big)+ \sigma_0\big(\sum_{i=1}^m L_{i}^*(v_{i,k}-v_{i,k+1})\big)\\
\operatorname{For \; \text{$i =1,\ldots,m$}}\\
\begin{array}{l}
\left\lfloor
\begin{array}{l}
w_{i,k+1} = \prox_{\sigma_i  g_{i}} \Big(w_{i,k}+\sigma_i O^{*}_iv_{i,k+1}\Big) + \sigma_iO^{*}_i(v_{i,k+1}-v_{i,k}).\\
\end{array}
\right.\\[2mm]
\end{array}\\
\end{array}
\right.\\[2mm]
\end{array}
\end{equation}
\end{theorem}
Then the following hold for some $(v_{1}^*,\ldots,v_{m}^*) \in \mathcal{D}_1$ and $x^*\in \mathcal{P}_1$. 
\begin{enumerate}
   \item \label{t:main2i}
   For every $j\in\{1,\ldots,m\}$, the sequence $(v_{j,k})_{k\in\NN}$ converges weakly to $v_{j}^*$.
    \item \label{t:main2ii} The sequence $(x_k)_{k\in\NN}$ converges weakly to $x^*$.
    \item\label{t:main2iii} If $f$ or $h$ is uniformly convex at $x^*$, then $x_k\to x^*$.
    \item\label{t:main2iv} If $g_j$ is uniformly convex at $v_{j}^*$, for some $j\in\{1,\ldots,m\}$, then $v_{j,k}\to v_{j}^*$.
    
\end{enumerate}

\begin{proof} Let us consider the problem of finding $\overline{\vv}\in\GG$ such that 
 \begin{equation}\label{e:rdu}
     0\in \nabla h_0(\overline{\vv}) +\sum_{i=0}^mK_{i}^*(\partial d_i\;\vuo\;\partial l_i)(K_i\overline{\vv}_i-r_{i}^*),
 \end{equation}
 which is an instance of the primal inclusion \eqref{e:pri}. Applying the proposed algorithm \eqref{algo:primal-dual} to \eqref{e:rdu} yields
 \begin{equation}
\label{algo:apprimal-dual}
\begin{array}{l}
\left\lfloor
\begin{array}{l}
\vv_{k+1} = \vv_k-\tau \big(\sum_{0=1}^mK_{i}^*w_{i,k}+ \nabla h_0 (\vv_k) -(r_1,\ldots,r_m)\big)\\
\operatorname{For \; \text{$i =0,\ldots,m$}}\\
\begin{array}{l}
\left\lfloor
\begin{array}{l}
w_{i,k+1} = J_{\sigma_i (\partial d_i)^{-1} } \Big(w_{i,k}+\sigma_i\big(K_i\vv_{k+1}-\nabla l_{i}^*(w_{i,k}) - r^{*}_i\big)\Big) + \sigma_i(K_i\vv_{k+1}-K_i\vv_k).\\
\end{array}
\right.\\[2mm]
\end{array}\\
\end{array}
\right.\\[2mm]
\end{array}
\end{equation}
We note that  $\nabla h_0\colon \vv\mapsto (W_i\nabla\ell^{*}_iW_{i}^*(v_{i})+r_i)_{1\leq i\leq m}$. Moreover, $\nabla l_{0}^* = \nabla h$ and, for every $i\in\{1,\ldots,m\}$, $\nabla l^{*}_i =0. $ Furthermore, 
for every $i\in\{1,\ldots,m\}$, $K_{i}^* = \pi_{i}^*\circ O_i\colon w_i\mapsto (0,\ldots,0,O_iw_i,0,\ldots 0)$, and $K_0=-\LL^*$ with $K_{0}^*=-\LL$.
Let $k\in\NN$, set $w_{0,k}=x_k$ and $\vv_k=(v_{1,k},\ldots, v_{m,k})$. Since $J_{\sigma_i (\partial d_{i})^{-1} } = \prox_{\sigma_i g_i}$, we observe that  \eqref{algo:apprimal-dual0} is exactly equivalent to \eqref{algo:apprimal-dual}. 

\ref{t:main2i}.
In view of Theorem \ref{t:main1}, we need to check the conditions ensuring the convergence of \eqref{algo:apprimal-dual0}. 
By Lemma \ref{l:lips}, the cocercivity constant of $C=\nabla h_0$ is $\beta_C$-cocoercive with
\begin{equation}
    \beta_C= \Bigg(\max_{1\leq i\leq m}\mu_{i}^{-1}\|W_i\|^2 \Bigg)^{-1}.
\end{equation}
Since $\nabla h$ is $\beta^{-1}$-Lipschitz, $\nabla l_0$ is  $\beta$-cocoercive. Therefore, \eqref{e:constep1} shows that \eqref{e:constep} is satisifed.
We first consider the condition \eqref{e:ran}. There exists $\overline{x}\in\HH$ such that 
\begin{equation}
    z\in \partial f(\overline{x}) + \sum_{i=1}^mL_{i}^*\big(\partial(O_i \triangleright g_i)\;\vuo\;\partial(W_i \triangleright\ell_i)\big)(L_i\overline{x}-r_i) + \nabla h(\overline{x}).
\end{equation}
This inclusion further implies the existence of $\overline{\vv} = (\overline{v}_1,\ldots,\overline{v}_m)\in\GG$ such that, for every $i\in\{1,\ldots,m\}$,
\begin{equation}
    \begin{cases}
        z-\sum_{i=1}^mL_{i}^*\overline{v}_i \in(\partial f+\nabla h)(\overline{x})\\
        \overline{v}_i \in \big(\partial(O_i \triangleright g_i)\;\vuo\;\partial(W_i \triangleright\ell_i)\big)(L_i\overline{x}-r_i)
    \end{cases}
    \Longleftrightarrow
    \begin{cases}
        \overline{x} \in (\partial f+\nabla h)^{-1}(z-\LL^*\overline{\vv})\\
        L_i\overline{x}-r_i\in \partial(O_i \triangleright g_i)^*(\overline{v_i} )+\partial(W_i \triangleright\ell_i)^*(\overline{v}_i),\notag
    \end{cases}
\end{equation}
which implies that, 
\begin{equation}
    0\in - \LL  (\partial f+\nabla h)^{-1}(z-\LL^*\overline{\vv})
    +\partial\Bigg(\bigoplus_{i=1}^m (O_i \triangleright g_i)^*\Bigg)(\overline{\vv})
    +\nabla\Bigg(\bigoplus_{i=1}^m (W_i \triangleright \ell_i)^*+\scal{\cdot}{r_i}\Bigg)(\overline{\vv}).\label{e:fermatd}
\end{equation}
By Fermat's rule, $\overline{\vv}$ solves the dual problem \eqref{e:dual1}. In view of the notation defined in Remark \ref{Re:1}, we can rewrite \eqref{e:fermatd} as 
\begin{equation}
    0 \in \ran\Bigg(\sum_{i=0}^m K_{i}^* (\partial d_{i}\;\vuo\;\partial l_i)(K_i\cdot-r_{i}^*)+\nabla h_0\Bigg),
\end{equation}
which shows that the condition \eqref{c:1} is satisfied. By Theorem \ref{t:main1iii}, the sequence $(\vv_k)_{k\in\NN}$ converges weakly to a point $\vv^*\in\GG$ that solves the inclusion 
\begin{equation}\label{e:pr2}
   0\in  \sum_{i=0}^m K_{i}^* (\partial d_{i}\;\vuo\;\partial l_i)(K_i\vv^*-r_{i}^*)+\nabla h_0(\vv^*).
\end{equation}
Therefore, $\vv^*\in\mathcal{D}_1$. 
\ref{t:main2ii}. By Theorem \ref{t:main1iii},
the sequence 
$(w_{0,k},\ldots,w_{m,k})_{k\in\NN}$ converges weakly to $(w_{0}^*,\ldots, w_{m}^*) \in\HH\times\YY_1\times \ldots,\YY_m$ that solve the dual of \eqref{e:pr2},
\begin{equation}
\label{e:dual2}
\big(\forall i\in\{0,\ldots,m\}\big) \; -r^{*}_i \in -K_i(\nabla h_0)^{-1}\bigg(-\sum_{i=0}^m K_{i}^* w^{*}_i\bigg) + (\partial d_{i})^{-1} w_{i}^* 
+(\partial l_{i})^{-1} w_{i}^*.
\end{equation}
Let us set $\uu\in (\nabla h_0)^{-1}\bigg(-\sum_{i=0}^m K_{i}^* w^{*}_i\bigg)$. Then, we have 
\begin{align}
    -\sum_{i=0}^m K_{i}^* w^{*}_i = \nabla h_0(\uu),
\end{align}
which is equivalent to 
\begin{equation}
    (\forall i\in\{1,\ldots,m\})\; L_i w_{0}^*- O_iw_{i}^*= \nabla (\ell_{i}^*\circ W^{*}_i)(u_i) +r_i.\label{e:Lie1}
\end{equation}
Note that, at $i=0$, we obtain 
\begin{equation}\label{e:Lie0}
    z \in \sum_{i=1}^m L^*_iu_i+ \partial f(w_{0}^*) + \nabla h(w_{0}^*),
\end{equation}
and for $i\in\{1,\ldots,m\}$, we obtain 
\begin{align}
    0\in - O_{i}^*u_i+ \partial g_{i}(w_{i}^*) &\Longleftrightarrow w_{i}^* \in \partial g_{i}^*(O_{i}^*u_i)\notag\\
    & \Longrightarrow O_iw_{i}^* \in O_i\partial g_{i}^*(O_{i}^*u_i )\notag\\
    &\Longleftrightarrow O_iw_{i}^* \in \partial (g_{i}^*\circ O_{i}^*)(u_i).\label{e:Lie2}
\end{align}
Adding \eqref{e:Lie1} and \eqref{e:Lie2}, we obtain
\begin{equation}
    (\forall i\in\{1,\ldots,m\})\; L_iw_{0}^* \in \partial (g_{i}^*\circ O_{i}^*)(u_i) + \nabla (\ell_{i}^*\circ W^{*}_i)(u_i) +r_i,
\end{equation}
which is equivalent to 
\begin{align}
     (\forall i\in\{1,\ldots,m\})\; u_i &\in \bigg( \partial (g_{i}^*\circ O_{i}^*) + \nabla (\ell_{i}^*\circ W^{*}_i)\bigg)^{-1}(L_iw_{0}^*-r_i)\notag\\
     &=\big(\partial(O_i \triangleright g_i)\;\vuo\;\partial(W_i \triangleright\ell_i)\big)(L_iw_{0}^*-r_i).\label{e:Lie3}
\end{align}
Combining \eqref{e:Lie3} with \eqref{e:Lie0} to obtain that $w_{0}^*$ is a solution to 
\begin{equation}
      z \in \sum_{i=1}^m L^*_i\big(\partial(O_i \triangleright g_i)\;\vuo\;\partial(W_i \triangleright\ell_i)\big)(L_iw_{0}^*-r_i)+ \partial f(w_{0}^*) + \nabla h(w_{0}^*). 
\end{equation}
Thus, by \cite[Proposition 16.6(ii)]{livre1} and Fermat's rule \cite[Theorem 16.3]{livre1}, $w_{0}^*$ is a solution to the primal problem \eqref{e:pri1}. Therefore, the sequence $(x_k)_{k\in\NN}$ 
converges weakly to $x^*= w_{0}^*$.

\ref{t:main2iii}\&\ref{t:main2iv}. The results follow directly from \ref{t:mainiv} \&\ref{t:mainv}.
\end{proof}
\begin{remark} The algorithm \eqref{algo:apprimal-dual0} is full splitting, and it does not require the proximity operators of the infimal postcompositions $O\triangleright g$ (we omit the index $i$ for simple), which are not tractable generally. However, in some special cases, it is possible to find a closed-form expression of $\prox_{O\triangleright g}$ by formula \cite{Luis21,Luis24},
\begin{equation}
    \prox_{O\triangleright g}\colon x\mapsto O\prox_{g,O}\;
    \text{where}\; 
    \prox_{g,O} \in 
 \underset{y \in \mathcal{Y}}{\operatorname{Argmin}}\big(f(y)+\frac{1}{2}\|Ox-y\|^2\big).
\end{equation}
Moreover, under the condition $ 0 \in \sri(\dom g^*-\ran O^*)$, we also have \cite{Luis21},  for any strictly positive $\sigma$,
\begin{equation}
    \prox_{\sigma (O\triangleright g)^*} = \prox_{\sigma g^*\circ O^*} = \Id-\sigma O\prox_{g,O}(\cdot/\sigma). 
\end{equation}
In these special cases, one apply directly Algorithm \ref{algo:primal-dual}. Precisely, let $\tau$ and $(\sigma_i)_{1\leq i\leq m}$ be strictly positive real numbers. Let $x_0\in\HH$ and, for each $i\in \{1,\ldots,m\}$, let $v_{i,0}$ in $\GG_i$. Iterate
\begin{equation}
\label{algo:primal-dual3a}
\begin{array}{l}
\left\lfloor
\begin{array}{l}
x_{k+1} = \prox_{\tau  f} \Big(x_k-\tau \big(\sum_{i=1}^mL_{i}^*v_{i,k}+ \nabla h(x_k) -z\big)\Big)\\
\operatorname{For \; \text{$i =1,\ldots,m$}}\\
\begin{array}{l}
\left\lfloor
\begin{array}{l}
q_{i,k} = v_{i,k}+\sigma_i\big(L_ix_{k+1}-W_i\nabla\ell_{i}^*(W^{*}_iv_{i,k}) - r_i\big)\\
v_{i,k+1} = q_{i,k}-O_i\sigma_i\prox_{\sigma_i g_i,O_i} (q_{i,k}/\sigma_i) + \sigma_i(L_ix_{k+1}-L_ix_k).\\
\end{array}
\right.\\[2mm]
\end{array}\\
\end{array}
\right.\\[2mm]
\end{array}
\end{equation}
\end{remark}
Suppose that the condition \eqref{e:ran} is satisfied and that 
\begin{equation}\label{e:constep3}
    \begin{cases}
     \big(\forall i\in\{1,\ldots,m\}\big)\;   \sigma_i < \mu_i/\|W_i\|^2,\\
        \frac{1}{2\tau}-\frac{1}{4\beta}-
\sum_{i=1}^m\sigma_i\|L_i\|^2\Big( \frac{\sigma_i}{2\mu_i}+\frac{2\mu_i}{\mu_i-\sigma_i}-\frac{3}{2}\Big) > 0.
    \end{cases}
\end{equation}
Then the following hold for some $(v_{1}^*,\ldots,v_{m}^*) \in \mathcal{D}_1$ and $x^*\in \mathcal{P}_1$. 
\begin{enumerate}
   \item \label{t:rm2i}
   For every $j\in\{1,\ldots,m\}$, the sequence $(v_{j,k})_{k\in\NN}$ converges weakly to $v_{j}^*$.
    \item \label{t:rm2ii} The sequence $(x_k)_{k\in\NN}$ converges weakly to $x^*$.
    \item\label{t:rm2iii} If $f$ or $h$ is uniformly convex at $x^*$, then $x_k\to x^*$.
    \item\label{t:rm2iv} If $g^{*}_j\circ O^{*}_j$ or $\ell^{*}_j\circ W^{*}_j$ is uniformly convex at $v_{j}^*$, for some $j\in\{1,\ldots,m\}$, then $v_{j,k}\to v_{j}^*$.
\end{enumerate}
The proof of these conclusions follows the same argument as presented in the proof of \cite[Theorem 4.2]{plc6}.
\begin{remark} Here are some connections to existing work.
\begin{enumerate}
\item Structured convex optimization problems involving infimal postcompositions implicitly appeared in \cite{Becker} where the infimal postcompositions used in term of their conjugatate functions. 
Its applications to image restoration problems can be found in \cite{Becker,Luis24}.

    \item When, for every $i\in\{1,\ldots,m\}$, $O_i\triangleright g_i = g_i$ and $W_i\triangleright \ell_i= \ell_i$, Problem \ref{prob:app} reduces to the problem considered in \cite{plc6}, which was subsequently investigated in \cite{vu2013,Condat,Bot2013,MinYan} and the references therein.
    \item The algorithm \eqref{algo:primal-dual3a} recovers the one presented in \cite[Section 4]{Tam19} as a special case where $h=0$, $m=1$, $\ell_1=0$ and $O_1=\Id$. In this special case, algorithm \eqref{algo:apprimal-dual0} reduces to an alternative form of the algorithm in \cite[Section 4]{Tam19}.
\end{enumerate}   
\end{remark} 
\begin{remark}
Our approach showed that each existing method solving Problem \ref{p:Prob0} with $Q=0$, can use to solve Problem \ref{prob:app}.
\end{remark}
\noindent{\bf Declaration.} The authors contributed equally. The authors report generative AI was not used in their research or preparation of this manuscript.

\end{document}